\documentstyle[12pt,amsfonts,epsfig]{article}

\newcommand{\new}{\newcommand}

\newcommand{\BOX}{\hbox {$\sqcap$ \kern -1em $\sqcup$}}
\newcommand{\qed}{\hskip 3em \hbox{\BOX} \vskip 2ex}
\new{\tnsr}{\otimes} \new{\iso}{\cong} \new{\union}{\cup}
\new{\implies}{\Rightarrow} \new{\maps}{\colon} \new{\A}{{\cal A}}
\new{\G}{{\cal G}} \new{\R}{{\Bbb R}} \new{\C}{{\Bbb C}}
\new{\tr}{{\rm tr}} \new{\range}{{\rm Range}} \new{\dom}{{\rm Domain}}
\new{\Lie}{{\rm Lie}} \new{\Fun}{{\rm Fun}} \new{\Aut}{{\rm Aut}}
\new{\U}{{\rm U}} \new{\abs}[1]{\left|#1 \right|}
\new{\norm}[1]{\left\| #1 \right\|} \new{\bracket}[1]{\langle #1
\rangle} \new{\defequals}{\stackrel{\rm def}{=}}
\new{\into}{\hookrightarrow} \new{\comb}[2]{{#1 \choose #2}}
\new{\lbl}[1]{\label{#1} \if\draft y
\smash{\makebox[0pt]{\hspace{-0.5in} \raisebox{8pt}{\rm\tiny #1}}} \fi
}

\new{\pic}[5]{\raisebox{#3pt}{
\hspace{#4pt}\epsfig{file=#1.eps,height=#2pt}\hspace{#5pt}}}

 \newtheorem{proposition}{Proposition}
\newtheorem{theorem}{Theorem}

\newtheorem{lemma}{Lemma} 
\newtheorem{corollary}{Corollary}

\newcounter{letter} \newcounter{numeral} \newcounter{Numeral}
\newenvironment{alphalist}{
\begin{list}{{\normalshape(\alph{letter})}}{\usecounter{letter}}
}{\end{list}}

\newenvironment{romanlist}{
\begin{list}{(\roman{numeral})}{\usecounter{numeral}}
}{\end{list}}

\begin{document}

      \begin{center}
      {\bf Functional Integration on Spaces of Connections \\}
      \vspace{0.5cm} {\em John C. Baez\\} \vspace{0.3cm} {\small
      Department of Mathematics \\ University of California\\
      Riverside, CA 92521-0135\\} \vspace{0.3cm}
      \renewcommand{\thefootnote}{\fnsymbol{footnote}} {\em Stephen
      Sawin\footnote{This author partially supported by NSF
      Postdoctoral Fellowship \#23068} \\} \vspace{0.3cm} {\small
      Department of Mathematics \\ Massachusetts Institute of
      Technology \\ Cambridge, MA 02139-4307\\ } \vspace{0.3 cm}
      {\small email: \tt baez@math.ucr.edu, sawin@math.mit.edu}
     \end{center}

\begin{abstract}

Let $G$ be a compact connected Lie group and $P \to M$ a smooth
principal $G$-bundle.  Let a `cylinder function' on the space $\A$ of
smooth connections on $P$ be a continuous function of the holonomies
of $A$ along finitely many piecewise smoothly immersed curves in $M$,
and let a generalized measure on $\A$ be a bounded linear functional
on cylinder functions.  We construct a generalized measure on the
space of connections that extends the uniform measure of Ashtekar,
Lewandowski and Baez to the smooth case, and prove it is invariant
under all automorphisms of $P$, not necessarily the identity on the
base space $M$.  Using `spin networks' we construct explicit functions
spanning the corresponding Hilbert space $L^2(\A/\G)$, where $\G$ is
the group of gauge transformations.
\end{abstract}

\section{Introduction}

Integrals over spaces of connections play an important role in modern
gauge theory, but as these spaces are infinite-dimensional, it is
often difficult to make heuristic computations involving such
integrals rigorous.  Suppose one has a smooth principal $G$-bundle $P
\to M$, and let $\A$ be the space of smooth connections on $M$.  Then
$\A$ is an affine space, and becomes a vector space after an arbitrary
choice of some point as origin, so initially it may be tempting to
integrate functions using some sort of `Lebesgue measure' on $\A$.
Unfortunately, various theorems \cite{BSZ} indicate that there are no
well-behaved translation-invariant measures on an infinite-dimensional
vector space.

One might then restrict ones ambition to integrating `cylinder
functions' and certain limits thereof.  A cylinder function on $\A$ is
one that depends on finitely many coordinates, that is, one of the
form
\[ F(A) = f(\ell_1(A), \dots, \ell_n(A)) \] where $\ell_i\maps \A \to
\R$ are continuous linear functionals and $f \maps \R^n \to \C$ is
bounded and continuous.  To integrate these all one needs is a
`cylinder measure'; the theory of these is well-developed and widely
used in probability theory, quantum mechanics and quantum field theory
\cite{BSZ,K}.

However, in gauge theory the fact that $\A$ is an affine space is in
many ways less important than the fact that the group $\G$ of gauge
transformations acts on it.  For example, one is often interested in
integrating gauge-invariant functions on $\A$.  This amounts to doing
integrals on $\A/\G$, an infinite-dimensional space which is not an
affine space, nor even a manifold, but a kind of stratified space.  In
applications to physics \cite{Baez4,Bruegmann,GP,L1} one is often
interested in integrating `Wilson loops', gauge-invariant functions of
the form
\[ F(A) = \tr({\rm T} \exp \int_\gamma A) \] where $\gamma$ is a
smooth loop in $M$, ${\rm T} \exp \int_\gamma A$ denotes the holonomy
of $A$ around $\gamma$, and the trace is taken in some
finite-dimensional representation of $G$.  Wilson loops are typically
not easy to approximate by cylinder functions unless $G$ is abelian,
so it is difficult to integrate them against cylinder measures.

Motivated by Rovelli and Smolin's work \cite{RS} on the loop
representation of quantum gravity, Ashtekar and Isham \cite{AI}
introduced a nonlinear version of the cylinder measure idea which is
specially adapted to gauge theory.  Taking advantage of subsequent
reworkings, we may describe their idea as follows. One first {\it
redefines} a `cylinder function' on $\A$ to be one of the form
\begin{equation}
  F(A) = f({\rm T} \exp \int_{c_1} A, \dots, {\rm T}\exp \int_{c_n} A)
\label{cyl}
\end{equation} where $c_i$ are piecewise smooth paths in $M$, the
holonomy ${\rm T} \exp \int_{c_i} A$ of the connection $A$ along $c_i$
is identified with an element of $G$ by means of a trivialization of
$P$ over the endpoints of $c_i$, and $f \maps G^n \to \C$ is
continuous.  Taking the completion of this algebra in the $\sup$ norm
\[ \|F\|_\infty = \sup_{A \in \A} |F(A)| ,\] one obtains a C*-algebra
of bounded continuous functions on $\A$.  Then one defines a
`generalized measure' $\mu$ on $\A$ to be a bounded linear functional
on this C*-algebra.  Using the Gelfand-Naimark spectral theory, this
C*-algebra can be identified with the C*-algebra of all continuous
functions on a compact space $\overline\A$ of which $\A$ is a dense
subset.  Elements of $\overline \A$ are called `generalized
connections', and the holonomy of one of these generalized connections
along a piecewise smooth path is still well-defined.  By the
Riesz-Markov theorem, generalized measures on $\A$ can be identified
with finite regular Borel measures on $\overline \A$.

One might hope for some generalized measure on $\A$ to serve as a
substitute for the nonexistent `Lebesgue measure' on $\A$.  At the
bare minimum one would like this generalized measure to be invariant
under all automorphisms of the bundle $P$ --- e.g., gauge
transformations and lifts of diffeomorphisms of $M$.  In a search for
something along these lines, Ashtekar and Lewandowski \cite{AL}
discovered that the study of generalized measures becomes more
manageable when one works with cylinder functions defined using
piecewise analytic paths.  Technically, the difficulty with piecewise
smooth paths is that they can intersect in very complicated ways ---
even in a Cantor set.  Piecewise analytic paths, on the other hand,
can only intersect in an infinite set if they overlap for some closed
interval.  This turns out to greatly simplify matters.

After further work by Ashtekar, Lewandowski and Baez
\cite{AL2,Baez,Baez2,Baez3,L}, the theory of generalized measures in
the analytic context now looks as follows.  One assumes $M$ is a
real-analytic manifold, $G$ is a compact connected Lie group, and $P
\to M$ is a smooth principal $G$-bundle.  One works only with cylinder
functions for which the paths $c_i$ are piecewise real-analytic.
Letting $\Fun_\omega(\A)$ denote the completion of this space of
cylinder functions in the $\sup$ norm, one then defines a `generalized
measure' to be a bounded linear functional on $\Fun_\omega(\A)$.

The results of Ashtekar and Isham still hold: $\Fun_\omega(\A)$ is
isomorphic to the C*-algebra of continuous functions on a space
$\overline \A$ containing $\A$ as a dense subset, and generalized
measures on $\A$ are the same as finite regular Borel measures on
$\overline \A$.  In the analytic context, however, it is not too hard
to construct a canonical generalized measure on $\A$, the `uniform'
generalized measure.  This generalized measure is invariant under all
automorphisms of the bundle $P$ that act on the base manifold $M$ as
real-analytic diffeomorphisms.  The uniform generalized measure is not
the only one invariant under all these automorphisms.  In fact, many
such generalized measures exist, and they can be constructed and ---
in a rather abstract sense --- classified using the notion of an
`embedded graph'.  An embedded graph $C$ is a finite set of analytic
paths $c_i \maps [0,1] \to M$ that are 1-1, embeddings when restricted
to $(0,1)$, and nonintersecting except possibly at their endpoints.
These paths are called the `edges' of the graph.  One can reduce the
study of holonomies along any finite set of real-analytic paths to the
case of a graph, because given any such set of paths, one can write
them as finite products of the edges of some graph (and their
inverses).  Given an embedded graph $C$ with $n$ edges, and
trivializing $P$ over the endpoints of all the edges, a generalized
measure $\mu$ on $\A$ determines a measure $\mu_C$ on $G^n$ by
\[ \int_{G^n} f(g_1, \dots, g_n) d\mu_C = \mu(F) \] where $F$ is
related to $f$ by equation (\ref{cyl}).  The measures $\mu_C$ for all
embedded graphs $C$ determine the generalized measure $\mu$, and the
uniform generalized measure on $\A$ is the unique one for which all
the measures $\mu_C$ are normalized Haar measure on some product of
copies of $G$.

It is natural to wonder whether these results depend crucially on the
use of analytic paths.  This is not a question of merely technical
interest.  One might argue that the analyticity assumptions are not so
bad, since every paracompact smooth manifold admits a real-analytic
structure, which is unique up to smooth diffeomorphism
\cite{Palais,Whitney}.  However, in applications to topological
quantum field theory and the loop representation of quantum gravity,
diffeomorphism-invariance plays a key role, and real-analytic
diffeomorphisms do behave very differently from smooth ones.  After
all, a real-analytic diffeomorphism of a connected manifold is
completely determined by its restriction to an arbitrarily small
neighborhood.  To see how this impinges on questions of real physical
interest, it is interesting to read the recent work of Ashtekar,
Lewandowski, Marolf, Mour\~ao and Thiemann on diffeomorphism-invariant
gauge theories \cite{ALMMT}.

The goal of this paper is to treat the case where $M$ is merely {\it
smooth}.  We work with cylinder functions on $\A$ for which the paths
$c_i$ are `curves', that is, {\em piecewise smoothly immersed}, and we
let $\Fun(\A)$ denote the completion of the space of these cylinder
functions in the $\sup$ norm.  For us, a generalized measure will be a
continuous linear functional on $\Fun(\A)$.  Note that if $M$ is
real-analytic then $\Fun_\omega(\A) \subseteq \Fun(\A)$, so any of our
generalized measures restricts to a generalized measure as defined in
the real-analytic context.  In particular, we construct a generalized
measure that is invariant under all automorphisms of the bundle $P$,
and which restricts to the uniform generalized measure when $M$ is
real-analytic.  We again call this the `uniform' generalized measure.

In fact, this uniform generalized measure was already constructed by
Ashtekar and Lewandowski \cite{AL} in the case $G = \U(1)$, using
special properties of abelian Lie groups.  From this point of view,
the advance of the present paper consists of being able to handle
nonabelian groups.  But our work also establishes a framework for
handling other generalized measures in the smooth context.

The main ideas behind this framework are as follows.  In analogy with
the analytic case, for every family of curves $C = \{c_1, \dots,
c_n\}$ a generalized measure $\mu$ on $\A$ determines a measure
$\mu_C$ on $G^n$ by
\[ \int_{G^n} f(g_1, \dots, g_n) d\mu_C = \mu(F) \] where $F$ is
related to $f$ by equation (\ref{cyl}).  The goal is thus to
reconstruct a generalized measure $\mu$ starting from such a measure
$\mu_C$ for every family $C$.  Of course, some conditions must hold
for such a collection of measures $\mu_C$ to come from a generalized
measure on $\A$.  In particular, not all $n$-tuples of elements of $G$
can be simultaneously attained as the holonomies of some fixed
connection along the curves in $C$, but only those lying in some
subset $\A_C \subseteq G^n$.  To come from a generalized measure, for
each family $C$ the measure $\mu_C$ will need to be supported on this
`attainable subset' $\A_C$, so we need a good understanding of this
subset.  In particular, in contrast to the analytic case, we cannot
reduce the problem to considering nice families such as embedded
graphs for which $\A_C = G^n$.

\begin{figure}[hbt]
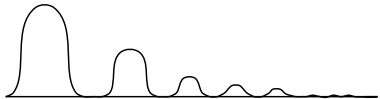

$$\pic{typeex}{35}{-15}{0}{0}$$
\caption{A family of curves with a type occurring infinitely often}
\label{fg:type}
\end{figure}

The reason why $\A_C$ may not be all of $G^n$ is that there may be
relations among the holonomies along the curves in the family $C$.
These relations occur when the curves overlap for some open interval,
so we need to introduce a notion of a `type' of possible overlap.  Due
to the complicated ways curves can intersect in the smooth context, a
given type may occur infinitely often in a family $C$; for a simple
example see Figure \ref{fg:type}.

The goal of Section 2 is to describe the possible holonomies of a
family of curves.  This is done first for especially well-behaved
families of curves called `tassels'.  Rough\-ly, a tassel based on a
point $p \in M$ is a family of curves for which, when it is restricted
to any neighborhood of $p$, the same types of overlap still occur.
This self-similarity forces $\A_C$ to be a subgroup of $G^n$, in fact
a subgroup easily presented in terms of the types of overlap occurring
in $C$.  Then we introduce the notion of a `web'.  This is a family
$W$ of curves that can be written as a union of tassels $W^1, \ldots,
W^k$, sufficiently separated so that $\A_W=\A_{W^1} \times \cdots
\times \A_{W^k}$.  We will show that the holonomies along any family
can be expressed in terms of the holonomies along some web, thus
giving an algebraic description of the possible holonomies.  The
proofs of these facts, that is of Propositions 1 and 2, can safely be
skipped by a reader looking for an initial overview of the results of
this paper.

In fact, if the family one started with consisted of the edges of an
embedded graph, the tassels this construction would produce would be
the restriction of the edges to each cell in a cell decomposition dual
to the graph.  Thus each tassel would contain one vertex $p$ of the
graph, and would in fact be based at $p$.  One should therefore think
of a web as a generalization of a finite graph, and a tassel based at
$p$ as a generalization of a neighborhood of a vertex $p$.

In Section 3 we give a criterion for a collection of measures $\mu_W$,
one for each web $W$, to arise from a generalized measure $\mu$ on
$\A$.  We also show that $\mu$ is uniquely determined by the measures
$\mu_W$, so that we have a tool for constructing generalized measures.
In Section 4 we apply this tool to construct the uniform generalized
measure.

In recent work on the loop representation of quantum gravity, `spin
network states' play an important role \cite{ALMMT,Baez,RS2,RS3}.
These have already been dealt with rigorously in the analytic context,
and in Section 5 we describe how they work in the smooth context.  The
basic idea is as follows.  Using the uniform generalized measure $\nu$
on $\A$, one can define a Hilbert space $L^2(\A)$ by completing
$\Fun(\A)$ in the norm associated to the inner product
\[ \langle F,G\rangle = \nu (\overline F G) .\] The group $\G$ of
gauge transformations acts on $\A$, and this gives rise to a unitary
representation of $\G$ on $L^2(\A)$.  We define $L^2(\A/\G)$ to be the
subspace of $\G$-invariant vectors in $L^2(\A)$.  The `spin network
states' form a very explicit `local' orthonormal basis of
$L^2(\A/\G)$, which is to say an orthonormal basis of the subspace
associated to each web $W$.  In the analytic context, they are formed
using embedded graphs whose edges are labeled with irreducible unitary
representations of $G$, and whose vertices are labeled with
intertwining operators from the tensor product of the representations
labeling the `incoming' edges, to the tensor product of the
representations labeling `outgoing' edges.  In the smooth context spin
networks are formed using webs equipped with similar, but more subtle,
representation-theoretic data.  An embedded graph is a special case of
a web, and in this case our spin network states reduce to the spin
network states as defined in the analytic context.  However, it is not
yet clear whether the spin networks can be combined in a simple
fashion to give an orthonormal basis of all of $L^2(\A/\G)$
simultaneously, as in the analytic case.

\section{Webs}

Fix a connected compact Lie group $G$, a smooth (paracompact)
$N$-dimensional manifold $M$, and a smooth principal $G$-bundle $P \to
M$.  By a {\em curve} we mean a piecewise smooth map from a finite
closed interval of $\R$ to $M$ that is an immersion on each piece.
Two curves are considered {\em equivalent} if one is the composition
of the other with a positive diffeomorphism between their domains.  A
{\em family} of curves is a finite set of curves with a chosen
ordering $c_1, \ldots ,c_n$.  If $C$ is such a family, let $\range(C)$
be the union of the ranges of the individual curves.

If $c_1\maps [a,b] \to M$ and $c_2\maps [c,d] \to M$ are two curves
such that $c_1(a)=c_2(d)$, we can form the {\em product\/} $c_1c_2$ by
gluing them together at this common point.  Of course this is defined
only up to equivalence of curves.  It is exactly like the product in
the fundamental groupoid, except that we do not identify homotopic
curves.  It is still associative, however, and there is a category
whose objects are points in $M$ and whose morphisms (other than
identity morphisms) are equivalence classes of curves.  Define the
{\em inverse\/} $c^{-1}$ of a curve $c$ to be $c$ reparametrized by an
order-reversing map, again up to equivalence.  This is not truly an
inverse for the product, but merely a contravariant functor.

If every curve in the family $C$ is equivalent to a (finite) product
of curves in the family $D$ and their inverses, we say that $C$ {\em
depends\/} on $D$.  We say that a collection of families of curves
$C^1, \dots, C^k$ is {\em independent\/} if when $i \ne j$, any curve
in the family $C^i$ intersects any curve in the family $C^j$, if at
all, only at their endpoints, and there is a neighborhood of each such
intersection point whose intersection with $\range(C^i \cup C^j)$ is
an embedded interval.  Obviously even if two families are not
independent, one may not depend on the other.

The above definitions are motivated by considering holonomies of
connections along these curves.  The map from curves to holonomies
given by such a connection sends product to product and inverse to
inverse.  If one family of curves depends on another, one can compute
the holonomy of a connection along all the curves in the first from
the same information about the second.  If two families are
independent, knowing the holonomies along one family tells one nothing
about the holonomies along the other.

A {\em subcurve\/} of a curve $c$ is a curve equivalent to the
restriction of $c$ to a subinterval of its domain.  The {\em
restriction\/} of a family $C$ to a closed set $K \subset M$ is the
family gotten by restricting each $c_i$ to each interval of
$c_i^{-1}[K]$.  A point $p \in \range(C)$ is a {\em regular point\/}
if it is not the image of an endpoint or nondifferentiable point of
$C$, and there is a neighborhood of it whose intersection with
$\range(C)$ is an embedded interval.

A family of curves $C$ is {\em parametrized consistently\/} if each
curve is parametrized so that $c_i(t)=c_j(s)$ implies $t=s$.  Thus
each of the curves is actually an embedding, and each point $p$ in the
range of the family is associated to a unique value of the parameter,
which we call $t(p)$.  If a family $\{c_1, \ldots, c_n\}$ is
parametrized consistently and $p$ is a point in $\range(C)$, define
the {\em type\/} of a regular point $p$, $\tau_p$, to be the Lie
subgroup of $G^n$ consisting of all $n$-tuples $(g_1, \dots, g_n)$
such that for some $g \in G$ we have $g_i = g$ if $p$ lies on $c_i$,
and $g_i = 1$ otherwise.  This gives a canonical isomorphism between
any type and $G$.

\begin{figure}[hbt]
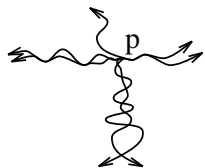

$$\pic{tassel}{60}{-35}{0}{0}$$
\caption{A tassel based at $p$}\label{fg:example}
\end{figure}

A fundamental concept in all that follows is that of a `tassel'.  A
family of curves $T$ is a {\em tassel based on\/} $p \in \range(T)$
if:
\begin{alphalist}
\item $\range(T)$ lies in a contractible open subset of $M$
\item $T$ can be consistently parametrized in such a way that $c_i(0)
= p$ is the left endpoint of every curve $c_i$
\item Two curves in $T$ that intersect at a point other than $p$
intersect at a point other than $p$ in every neighborhood of $p$
\item Any type which occurs at some point in $\range(T)$ occurs in
every neighborhood of $p$.
\end{alphalist} One may visualize the curves of the tassel as
radiating outwards from the base $p$.  See Figure \ref{fg:example} for
an example.

Finally, a {\em web} $W$ is a finite independent collection of tassels
$W^1, \dots, W^k$ such that no tassel contains the base of another.
We frequently apply concepts defined for families of curves to webs
without comment, using the fact that the web $W$ has an associated
family $W^1 \cup \cdots \cup W^k$.  For example, we say that a family
depends on a web $W$ if it depends on the family $W^1 \cup \cdots \cup
W^k$.  Our first main result about webs is:

\begin{proposition}\label{web1}  Any family of curves $C$ depends on a
  web $W$.
\end{proposition}

We begin with a technical lemma.

\begin{lemma} \label{lm:tech}
Let $C$ be a family of smooth curves $c_1, \ldots,c_n$.
\begin{alphalist}
\item The preimage of any point in $M$ under any $c_i$ is finite.

\item Every point $p \in \range(C)$ has a contractible open
  neighborhood $O$ admitting coordinates $x_1, \ldots x_N$ such that
  for each $i$, $dx_1(c_i(t))/dt \ne 0$ on $c_i^{-1}[O]$.

\item Given any point $p \in \range(C)$ and any open neighborhood $U$
  of $p$, there is an open subneighborhood $N$ of $p$ such that for
  each $i$, $c_i^{-1}[N]$ is a finite union of intervals, each
  containing a point of $c_i^{-1}[p]$.

\item The set of regular points is open and dense in $\range(C)$.

\item Given any point $p \in \range(C)$ and any open neighborhood $U$
  of $p$, there is an open subneighborhood $N$ with the properties in
  part (c) such that every point of $\range(C)$ lying on the boundary
  of $N$ is a regular point.
\end{alphalist}
\end{lemma}

{\it Proof.}

\begin{alphalist}
\item If not, the preimage would have an accumulation point, and at
  that point $c_i$ would not be an immersion.

\item We can choose an open neighborhood $U$ about $p$ with
  coordinates $x_1, \dots, x_N$ such that for all $i$ we have $d
  x_1(c_i(t))/dt \neq 0$ at all of the finitely many points in the
  preimage of $p$ under $c_i$.  Each such point has an open interval
  around it such that $dx_1(c_i(t))/dt \neq 0$ throughout that
  interval.  The union of the images of the complements of these
  intervals is a compact set $K \subseteq M$ not containing $p$.  It
  follows that any contractible open neighborhood $O$ of $p$ contained
  in $U - K$ has the desired properties.

\item Choose a coordinate patch $O$ around $p$ as in part (b) of this
  lemma, and consider the hyperplane through $p$ on which $x_1$ is
  constant.  The points of intersection of $\range(C)$ with this
  hyperplane are all transverse, so they have no accumulation points.
  Thus a small open neighborhood of $p$ in the hyperplane only
  intersects $\range(C)$ at $p$.  Shrinking this neighborhood to a
  sufficiently small subneighborhood, its product with a sufficiently
  small open interval in the $x_1$ axis is an open neighborhood $N$ of
  $p$ that only intersects each $c_i$ in finitely many embedded open
  intervals containing $p$.  This choice of $N$ has the desired
  properties.  See Figure \ref{fg:nbhd} for an illustration.

\begin{figure}[hbt]
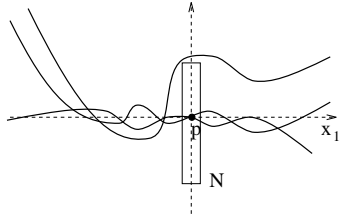

$$\pic{neighborhood}{80}{-35}{0}{0}$$
\caption{Choosing the neighborhood $N$ of $p$} \label{fg:nbhd}
\end{figure}

\item Consider any neighborhood $U$ of a point $p \in \range(C)$.  For
  each point in $U$, consider the total number of points in the
  preimages of all the curves $c_i$.  By part (a), we can pick a point
  $p_0 \in \range(C)$ for which this number is minimal.  We will show
  $p_0$ is regular.

  Choose a subneighborhood $N$ of $U$ as in part (c) of this lemma,
  and small enough that each $c_i$ is 1-1 on each component of the
  preimage of $N$.  Each $r \in N \cap \range(C)$ has at most one
  preimage point in each of the intervals comprising $c_i^{-1}[N]$,
  and since $p_0$ was minimal and has exactly one preimage point in
  each of these intervals by part (c), $r$ must have exactly one in
  each.  Thus the images of these intervals in $N$ must coincide, and
  hence $N \cap \range(C)$ must be an embedded open interval in $M$.
  Since $p_0$ cannot be the image of an endpoint, $p_0$ is regular.

\item Choose $N$ as in the proof of part (c).  Recall that $d
  x_1(c_i(t))/dt \neq 0$ for all $i$ and all $t$ in the preimage of
  $p$ under $c_i$.  If we choose the interval $(a,b)$ in the $x_1$
  axis used to define $N$ sufficiently small, each curve $c_i$
  intersects the boundary of $N$ only at the planes where $x_1$ equals
  $a$ or $b$.  Moreover, choosing this interval sufficiently small
  guarantees that the intersection points are transversal.  By part
  (d), it follows that we can choose $a$ and $b$ such that these
  intersection points are all regular points.  \qed

\end{alphalist}

{\em Proof of Proposition \ref{web1}.}  Let $C$ be a family of curves;
we may assume all the curves in $C$ are smooth, since any family
depends on a family of smooth curves.  By Lemma \ref{lm:tech}(b) and
the compactness of $\range(C)$, we can cover $\range(C)$ with open
sets $O_1, \ldots ,O_m$ such that each $O_l$ is contractible and
admits coordinates as in the lemma.

We claim that each $p \in \range(C)$ has an open neighborhood $N_p$
with the following properties:
\begin{romanlist}
\item $N_p$ is contained in every $O_l$ containing $p$.
\item $c_i^{-1}[\overline N_p]$ is a finite union of intervals.
\item $C$ restricted to $\overline N_p$ depends on a tassel based at
$p$.
\item the points of $\range(C)$ lying on the boundary of $N_p$ are all
regular points.
\end{romanlist}

To see this, we first tentatively take as $N_p$ the subneighborhood
given by Lemma \ref{lm:tech}(e) of the intersection of the $O_l$'s
containing $p$.  Then $c_i^{-1}[\overline N_p]$ is a finite union of
closed intervals, and the points of $\range(C)$ lying on the boundary
of $N_p$ are all regular points.  Use as coordinates on $N_p$ any of
the coordinates on the $O_l$'s containing $p$; without loss of
generality we assume that $x_1(p) = 0$.  By Lemma \ref{lm:tech}(b) the
restriction of $C$ to $\overline N_p$ is a family consistently
parametrized by the coordinate $x_1$, and by Lemma \ref{lm:tech}(e)
each curve $c$ in this family has $0$ in its domain and $c(0) = p$.
Take each curve in this family, and if its domain $[a,b]$ contains $0$
in its interior, replace it with the two subcurves given by its
restriction to $[0,b]$ and the inverse of its restriction to $[a,0]$.
Denote the resulting family by $C_p$.  We can parametrize each curve
in $C_p$ by $|x_1|$, and then not only will $C_p$ be consistently
parametrized, but also each curve $c$ in $C_p$ will have $c(0) = p$ as
its left endpoint.

Note that the family $C$ restricted to $\overline N_p$ depends on the
family $C_p$.  Thus $N_p$ has all the properties claimed except that
$C_p$ might not be a tassel.  By the previous paragraph, and since
$\range(C_p)$ lies in some contractible open set $O_l$, the only way
$C_p$ can fail to be a tassel is by violating conditions (c) or (d) in
the definition of a tassel.

\begin{figure}[hbt]
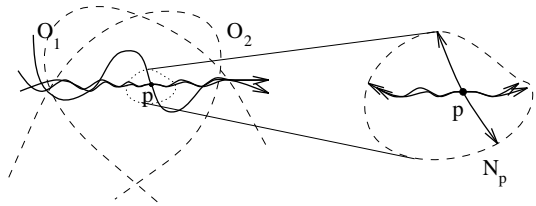

$$\pic{trim}{75}{-35}{0}{0}$$
\caption{Forming a tassel in a neighborhood of $p$}\label{fg:maketassel}
\end{figure}

To get condition (c) to hold, choose a neighborhood of $p$ in $N_p$
small enough that any two curves which intersect do so arbitrarily
close to $p$, and choose a subneighborhood as in Lemma
\ref{lm:tech}(e).  Use this subneighborhood as a new choice of $N_p$,
and restrict $C_p$ to the new $\overline N_p$.  This leaves us with a
neighborhood $N_p$ with all the properties claimed except that $C_p$
might violate condition (d) in the definition of a tassel.  To get
condition (d) to hold, note that for each type $\tau$ occurring in
$C_p$, either $\tau$ occurs at a sequence of points approaching $p$,
or it does not.  If it does not, choose a neighborhood of $p$ in $N_p$
which excludes all points of type $\tau$, and choose a subneighborhood
as in Lemma \ref{lm:tech}(e).  Use this subneighborhood as a new
choice of $N_p$, and restrict the original family $C_p$ to the new
$\overline N_p$, obtaining a new family $C_p$.  Since there were only
finitely many distinct types in the original $C_p$, and the new $C_p$
has fewer types, this process must end.  Now $C_p$ is a tassel and
$N_p$ is a neighborhood of $p$ having all the properties claimed.  See
Figure \ref{fg:maketassel} for an example of this construction.

\begin{figure}[hbt]
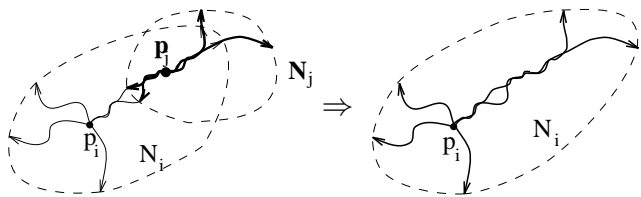

$$\pic{intersect}{70}{-30}{0}{0} \implies
\pic{union}{70}{-30}{0}{0} $$
\caption{The union is a tassel based at $p_i$}\label{fg:union}
\end{figure}

Now cover $\range(C)$ with finitely many of these open sets $N_p$.
Call them $N_1, \ldots,$ $N_k$, call the associated tassels $W^1,
\ldots, W^k$, and call the points at which they are based $p_1, \ldots
p_k$.  We claim that if $p_j \in N_i$ for some $i \neq j$, then $N_i
\cup N_j$ is still a neighborhood of $p_i$ with properties (i-iv).
Properties (ii) and (iv) are clear.  For (i), note that $N_i \cup N_j$
is contained in any $O_l$ containing $p_i$, because $N_i$ is and $p_j
\in N_i \subseteq O_l$ so $N_j$ is as well.  For (iii), coordinatize
$N_i \cup N_j$ using the coordinates on some $O_l$ containing $p_i$,
and construct a family as before, breaking the restriction of each
curve in $C$ to $\overline N_i \cup \overline N_j$ into two subcurves
with $x_1 \ge 0$ and $x_1 \le 0$ if necessary, and parametrizing them
consistently by the value of $|x_1|$.  To see that this family is a
tassel, the only nontrivial thing to check is condition (d).  Notice
that any type occurring in this family in $\overline{N}_i$ corresponds
to a type occurring in $W^i$, and therefore occurs arbitrarily close
to $p_i$.  Any type occurring in the family in
$\overline{N}_j-\overline{N}_i$ corresponds to a type occurring in
$W^j$, and thus occurs arbitrarily close to $p_j$.  But then it also
occurs in $N_i$, and thus arbitrarily close to $p_i$.  See Figure
\ref{fg:union} for an illustration.  Here bold curves are in $W^j$,
light curves are in $W^i$, and medium weight curves are in their
union.

Using this fact, we can replace the $N_i$ by unions thereof until no
$N_i$ contains $p_j$ for $i \neq j$, and succeed in covering
$\range(C)$ with open neighborhoods $N_1, \ldots , N_k$ of points
$p_1,\dots, p_k$ such that (i-iv) hold and such that $p_i \in N_j$
only for $i = j$.

To prove the theorem, we will now shrink each $N_i$ to a smaller
neighborhood of $p_i$ while maintaining properties (i-iv), so that the
resulting neighborhoods no longer intersect, but their closures still
cover $\range(C)$.  When we have done this, the restriction of $C$ to
the closure of each neighborhood will depend on a tassel by (iii).
Moreover, these tassels will form a web by (iv).  Finally, $C$ will
depend on this web by (ii).

To shrink the $N_i$ in this way, first replace each $N_j$ for $j >1$
with $N_j-\overline{N}_1$, leaving $N_1$ the same.  Properties (i) and
(iv) for the original $N_j$'s easily imply those properties for the
new $N_j$'s.  Property (ii) holds because the new $N_j$ has only
finitely many boundary points.  As for property (iii), recall that
$p_i \in N_j$ only for $i = j$.  the only way $W^j$ restricted to the
new $\overline N_j$ could fail to be a tassel is by having a component
that does not pass through $p_j$; but this could only happen if a
curve of $W^1$ had been a subcurve of $W^j$, in which case $p_1 \in
N_j$.

Next replace $N_j$ with $N_j-\overline{N_2}$ for each $j > 2$, and so
on.  When we are done, we find that $C$ depends on the modified $W^1,
\ldots, W^k$, which are all tassels, contain neighborhoods of the
$p_i$'s, and intersect only at boundary points.  Since all boundary
points are boundary points of the original $N_j$'s, they are regular
points of $C$, and therefore satisfy the condition for boundary points
of a web.  We thus obtain tassels $W^1, \dots, W^k$ forming a web on
which $C$ depends.  \qed

Let $\A$ be the space of smooth connections on $P$, equipped with its
$C^\infty$ topology.  Given a curve $c\maps [a,b] \to M$, let $\A_c$
be the set of functions
\[ \theta \maps P_{c(a)} \to P_{c(b)} \] compatible with the right
action of $G$ on $P$:
\[ \theta(xg) = \theta(x)g .\] Given $A \in \A$, the holonomy ${\rm T}
\exp \int_c A$ of $A$ along $c$ is an element of this set $\A_c$.  Of
course, a trivialization of $P$ at the endpoints of $c$ allows us to
identify $\A_c$ with $G$, and this gives $\A_c$ the structure of a
smooth manifold in a manner independent of the trivialization.  Note
also that $\A_c$ and the holonomy ${\rm T} \exp\int_c A \in \A_c$ only
depend on the equivalence class of $c$.

More generally, if $C = \{c_i \colon 1 \le i \le n\}$ is a family, let
\[ \A_C \subseteq \prod_{i=1}^n \A_{c_i} \] be the subspace consisting
of all elements of the form
\[ ({\rm T}\exp\int_{c_1} A, \dots ,{\rm T}\exp \int_{c_n} A) \] for
some connection $A \in \A$.  We call $\A_C$ the space of {\em
connections on} $C$, and give it the subspace topology.  If $W$ is a
web consisting of tassels $W^1, \dots, W^k$, we define $\A_W$ to be
$\A_{W^1 \cup \cdots \cup W^k}$, and again call this the space of
connections on $W$.

Note that the map $p_C \maps \A \to \A_C $ given by
\[ p_C(A) = ({\rm T}\exp\int_{c_1} A, \dots ,{\rm T}\exp \int_{c_n} A)
\] is continuous and onto.  Furthermore, if the product of curves
$c_1c_2$ exists there is a smooth map
\begin{eqnarray}
     \A_{c_1} \times \A_{c_2} &\to& \A_{c_1c_2} \nonumber \\
(\theta_1,\theta_2) &\mapsto& \theta_1\theta_2
.\nonumber\end{eqnarray} There is also for any curve $c$ a smooth map
\begin{eqnarray}         \A_c &\to& \A_{c^{-1}}  \nonumber \\
\theta &\mapsto& \theta^{-1}. \nonumber \end{eqnarray} Thus if $C$
depends on $D$, a particular choice of a way to write each curve in
$C$ as a product of curves in $D$ and their inverses gives a smooth
map
\[ p_{CD} \maps \A_D \to \A_C . \] Note that
\[ p_{CD} p_D = p_C. \] Since $p_D$ is onto, it follows that $p_{CD}$
is independent of how we write curves in $C$ in terms of curves in
$D$.  Since $p_C$ is onto, it also follows that $p_{CD}$ is onto.

Now suppose $\mu$ is a finite Borel measure on $\A$.  Then for any
family $C$, $\mu$ pushes forward by the map $p_C$ to a finite Borel
measure $\mu_C$ on $\A_C$.  The collection $\{\mu_C\}$ satisfies an
obvious consistency condition: whenever $C$ depends on $D$, the
measure $\mu_D$ pushes forward by the map $p_{CD}$ to the measure
$\mu_C$.  Not all collections of measures $\{\mu_C\}$ satisfying this
consistency condition arise from finite Borel measures on $\A$ in this
way, but as we shall see in Section 3, if such a collection satisfies
a certain uniform bound, it arises from a generalized measure on $\A$.
This is essentially how we construct the uniform generalized measure.
However, it is easier in practice to construct generalized measures
from collections $\{\mu_W\}$ where $W$ ranges over all webs, rather
than all families.  Proposition 1 is one of the results we need for
this, since it allows us to express any family in terms of a web.  The
second result we need is a good description of $\A_W$ when $W$ is a
web.

If $T$ is a tassel consisting of the curves $\{c_i\colon\; 1 \le i \le
n\}$, then we let $G_T$ be the smallest closed subgroup of $G^n$
containing all the types occurring in $T$, which of course is a Lie
group.  Then we have:

\begin{proposition}\label{web2}
If $T = \{c_i\colon\; 1 \le i \le n\}$ is a tassel, and we fix a
trivialization of $P$ over the endpoints of the curves $c_i$ to
identify $\A_T$ with a subset of $G^n$, then we have $\A_T = G_T$.  If
$W$ is a web consisting of tassels $W^1, \dots, W^k$, then $\A_W =
\A_{W^1} \times \cdots \times \A_{W^k}$.
\end{proposition}

{\em Proof.}  First suppose $T$ is a tassel.  Since $\range(T)$ is
contained in a contractible $U$, we can trivialize $P$ over
$\range(T)$, and by a suitable gauge transformation we choose this
trivialization so that it agrees with the specified trivialization of
$P$ over the endpoints of the curves $c_i$.  This allows us to treat
the holonomy of a connection along any of these curves from any point
$c_i(s)$ to any point $c_i(t)$ as an element of $G$.  It also allows
us to treat a connection on $P|U$ as a $\Lie(G)$-valued one-form.

We claim that given finitely many disjoint neighborhoods $N_\alpha
\subseteq U$ intersecting $\range(T)$ in open intervals $I_\alpha$
which contain no endpoints or nondifferentiable points of curves in
$T$, there is a connection $A_0$ on $P|U$ whose holonomy along
$I_\alpha$ is $g_\alpha$ for any $g_\alpha \in G$.  To see this, map
$I_\alpha$ to $G$ smoothly so that a neighborhood of its left endpoint
gets sent to $1$ and a neighborhood of its right endpoint gets sent to
$g_\alpha$.  Pull the derivative back to $I_\alpha$, extend it to a
smooth $\Lie(G)$-valued 1-form on $N_\alpha$ which is trivial in a
neighborhood $O_\alpha$ of the two endpoints, and multiply it by a
smooth function which is $1$ on $I_\alpha$ and $0$ near the boundary
of $N_\alpha$ outside of $O_\alpha$.  This gives a connection on
$P|N_\alpha$ whose holonomy along $I_\alpha$ is $g_\alpha$.  Defining
$A_0$ this way on each $N_\alpha$ and setting $A_0 = 0$ outside $N =
\bigcup N_\alpha$ proves the claim.

Notice if one of these intervals $I_\alpha$ is of type $\tau$,
$t[I_\alpha]\subset [a,b]$, and $[a,b]$ is in the domain of every
$c_i$, then the sequence of holonomies along the $c_i$ restricted to
$[a,b]$ can be made to be any element of $\tau$ by the above
procedure.

Consider any element of $G_T$ and write it as $\prod_{i=1}^n g_i$
where each $g_i$ is in $\tau_i$, a type occurring in $\range(T)$.  By
the definition of tassel, we can choose a decreasing sequence (with
respect to the parameter $t$) of regular points $p_i$ of type
$\tau_i$, for $i=1, \ldots, n$, such that each $t(p_i)$ is in the
interior of the domain of every curve in $C$.  Choose nonintersecting
neighborhoods $N_i$, and construct a connection $A_0$ on $P|U$ which
is trivial outside the $N_i$ and with holonomy $\prod_{i=1}^n g_i$.
Thus every element of $G_T$ is the holonomy of some connection $A_0$
on $P|U$.  Moreover, since $\range(T)$ is closed, we can find a
connection $A \in \A$ on all of $P$ which equals $A_0$ on $\range(T)$,
and thus has the same holonomy along each curve $c_i$.  It follows
that $G_T \subseteq \A_T$.

On the other hand, consider the map ${\cal C}:\R^+ \to \bigoplus_n TM$
sending each $t$ to $\bigoplus_n (c_i(t), c'_i(t))$.  If $c_i(t)$ is
not defined, use $(q_i,0)$, where $q_i$ is the right endpoint of
$c_i$, and if $c'_i(t)$ is not defined, use $(c_i(t),0)$.  This is
continuous except at finitely many points, namely endpoints or points
of nondifferentiability of any $c_i$.  Since $A$ gives a
$\Lie(G)$-valued one-form, we can interpret it as a map $A:\bigoplus_n
TM \to \bigoplus_n \Lie(G)$, so that $A \circ {\cal C}: \R^+
\to\bigoplus_n \Lie(G)$ is continuous except at finitely many points.

The set of $t$ such that $q$ is regular for all $q$ with $t(q)=t$ is
open dense, by Lemma \ref{lm:tech}(d).  For such a $t$, $ A\circ {\cal
C}$ is a sum of elements in the Lie algebras of the types occurring
with parameter value $t$, and thus is in $\Lie(G_T)$. By continuity
the range of $ A\circ {\cal C}$ is in $\Lie(G_T)$ except for finitely
many points.  But $p_T(A)$ is the endpoint of a curve in $G^n$
starting at the identity and having derivative $A \circ {\cal C}(t)$
at $t$.  Since this curve lies entirely in $G_T$, $p_T(\A) \in G_T$,
so $\A_T \subseteq G_T$.

Now suppose $W$ is a web consisting of tassels $W^1, \dots, W^k$.
Clearly $\A_W \subseteq \A_{W^1} \times \A_{W^k}$, so we need merely
prove the reverse inclusion.  In fact, we shall fix a trivialization
of $P$ over the subset of $M$ consisting of all the endpoints of the
curves in $W$, so as to identify each space $\A_{W^j}$ with $G_{W^j}$,
and we shall show that given $(g_1, \cdots, g_k) \in G_{W^1} \times
\cdots G_{W^k}$, there is a connection $A \in \A$ with $p_{W ^j}(A) =
g_j$ for all $i$.  Let $U_j$ be a contractible neighborhood containing
the tassel $W^j$. We can choose a trivialization of $P$ over each
$U_j$ which agrees with the above trivialization over all the
endpoints of the curves in $W$.  Moreover, we can choose these
trivializations so that they agree in a small neighborhood $O$ of all
these endpoints.  The construction above then gives for each tassel
$W^j$ a connection $A_j$ on $P|U_j$ such that $A_j$ has the desired
holonomies along all curves in $W^j$. Moreover, given any neighborhood
$V_j$ of $\range(W^j)$ with $\overline V_j \subset U_j$, we can choose
$A_j$ such that $A_j$ vanishes in $V_j$ except in an arbitrarily small
neighborhood $N_j$ of finitely many points in the interiors of the
curves in $W^j$.  (Here we use the trivialization of $P|U_j$ to think
of $A_j$ as a $\Lie(G)$-valued 1-form.)  If we choose the $V_j$'s
small enough that $V_i \cap V_j \subseteq O$ for $i\ne j$, and choose
the $N_j$'s small enough that $V_i \cap V_j \cap N_k = \emptyset$ for
all $i\ne j$, then the connections $A_j$ agree on all the overlaps
$V_i \cap V_j$ so there exists a connection $A_0$ on $P|\bigcup V_i$
having the desired holonomies on all curves in every tassel $W^j$.
Since $\range(W) \subset \bigcup V_j$ is closed there exists a
connection $A \in \A$ that equals $A_0$ over $\range(W)$, so
$p_{W^j}(A) = g_j$ for all $j$ as desired.  \qed

\section{Generalized Measures}

Let $\Fun_0(\A)$ be the algebra of {\it cylinder functions} on $\A$,
that is functions of the form
\[ F(A) = f(p_C(A)) \] where $C$ is some family of curves and $f \maps
\A_C \to \C$ is continuous.  Let $\Fun(\A)$ be the completion of
$\Fun_0(\A)$ in the $\sup$ norm.  We define a {\em generalized
measure} on $\A$ to be a continuous linear functional on $\Fun(\A)$.
Given a generalized measure $\mu$ on $\A$, for any family (or web) $C$
we can define a bounded linear functional $(p_C)_{*}\mu$ on the
algebra of continuous functions on $\A_C$ by:
\[ ((p_C)_* \mu)(f) = \mu(f \circ p_C) .\] By the Riesz-Markov
theorem, such a bounded linear functional is just a finite regular
Borel measure on $\A$.  (Henceforth when we write simply `measure' we
shall always mean a finite regular Borel measure.)

In short, a generalized measure on $\A$ determines a collection of
measures on the spaces $\A_C$ for all families $C$, and in fact such a
collection satisfying certain conditions uniquely determines a
generalized measure.  In light of Propositions \ref{web1} and
\ref{web2}, however, it is natural to translate this into the language
of webs.

\begin{theorem} \label{thm:gen}
Given a generalized measure $\mu$ on $\A$ and setting $\mu_W = (p_W)_*
\mu$ for any web $W$, the collection $\{\mu_W\}$ is:

\begin{alphalist}

\item Consistent: if the web $W$ depends on the web $X$ then
$(p_{WX})_* \mu_X = \mu_W$.

\item Uniformly bounded: the linear functionals $\mu_W \maps C(\A_W)
\to \C$ are uniformly bounded as $W$ ranges over all webs.

\end{alphalist}

Conversely, given any such consistent and uniformly bounded collection
$\{\mu_W\}$ of measures on the spaces $\A_W$, there exists a unique
generalized measure $\mu$ on $\A$ for which $(p_W)_* \mu = \mu_W$ for
all webs $W$.
\end{theorem} {\em Proof.} It is clear that given a generalized
measure $\mu$ on $\A$, the measures $\mu_W = (p_W)_* \mu$ are
consistent, and are uniformly bounded by the norm of $\mu$.

For the converse, suppose $\{\mu_W\}$ is a collection of measures on
the space $\A_W$ satisfying (a) and (b).  First we define a linear
functional $\mu_0$ on $\Fun_0(\A)$ as follows.  For any $F \in
\Fun_0(\A)$, choose a family $C$ and let $f_C \maps \A_C \to \C$ be a
continuous function with $F = f_C\circ p_C$.  By Proposition
\ref{web1}, there is a web $W$ upon which $C$ depends, so defining $f=
f_C \circ p_{CW}$, we have $F=f \circ p_W$. Now define
\[ \mu_0(F) = \mu_W(f).  \] Of course, we need to check that $\mu_0$
is well-defined and linear.  Suppose that $f' \maps \A_{W'} \to \C$ is
also continuous and $F = f' \circ p_{W'}$.  By Proposition \ref{web1}
again choose $X$ upon which $W \cup W'$, and hence $W$ and $W'$,
depend.  Then by (a)
\[ \mu_W = (p_{WX})_* \mu_X, \qquad \mu_{W'} = (p_{W'W})_* \mu_X.\]
 Also, since $f \circ p_W = f' \circ p_{W'}$ we have $f \circ p_{WX}
 \circ p_X = f \circ p_{W'X} \circ p_W$, but since $p_X$ is onto this
 implies
\[ f \circ p_{WX} =f' \circ p_{W'X} .\] Thus we have
\begin{eqnarray}   \mu_W(f) &=& ((p_{WX})_* \mu_X)(W)\nonumber\\
&=& \mu_X(f \circ p_{WX}) \nonumber\\ &=& \mu_{X}(f' \circ p_{W'X})
\nonumber\\ &=& ((p_{W'X})_* \mu_X)(f')\nonumber\\ &=& \mu_{W'} (f')
\nonumber \end{eqnarray} so $\mu_0$ is well-defined.  The linearity of
$\mu_0$ then follows from the linearity of each of the $\mu_W$'s.

By (b) we can choose $M > 0$ such that $\|\mu_W\| < M$ for all $W$,
and this implies that $|\mu_0(F)|$ for all $F \in \Fun_0(\A)$.  Since
$\Fun_0(\A)$ is dense in $\Fun(\A)$, $\mu_0$ extends uniquely to a
bounded linear functional $\mu$ on $\Fun(\A)$.  By construction,
$(p_W)_* \mu = \mu_W$ for all $W$.  The uniqueness of $\mu$ with this
property also follows from the fact that $\Fun_0(\A)$ is dense in
$\Fun(\A)$.  \qed

In fact, generalized measures on $\A$ are the same thing as measures
on the projective limit $\overline \A$ of the spaces $\A_C$, where the
families $C$ are ordered by dependence.  In these terms, Proposition
\ref{web1} says that webs are cofinal in the set of all families, and
Theorem \ref{thm:gen} is seen as a special case of a very general
result, namely that a measure on a projective limit of spaces can be
constructed from a consistent and uniformly bounded collection of
measures on any cofinal set of these spaces.  Ashtekar and Lewandowski
have given a clear exposition of this approach in the analytic context
\cite{AL2}, but here we chose to prove everything `from scratch.'

Elements of $\overline \A$ may be called {\em generalized connections}
on $P$.  Abstractly, $\overline \A$ is simply the Gelfand spectrum of
the C*-algebra $\Fun(\A)$.  The space $\A$ of smooth connections on
$P$ naturally maps into $\overline \A$ in a one-to-one and continuous
way, and the image is dense in $\overline \A$.  Thus generalized
connections may be regarded as limits of smooth connections.

\section{The Uniform Measure}

In this section we construct a generalized measure $\nu$ on $\A$ which
we call the `uniform measure'.  Theorem \ref{thm:gen} suggests that we
do this by choosing for each web $W$ a measure $\nu_W$ in some
canonical way.  In the special case of a web consisting of a single
tassel $T$, fixing a trivialization over the endpoints lets us think
of $\nu_T$ as a measure on $G_T$.  Since $G_T$ is a compact Lie group,
an obvious choice is Haar measure on $G_T$.  For more general webs it
is natural to use a product of Haar measures.  This in fact gives a
generalized measure.

\begin{theorem}\label{thm:uni} There exists a unique generalized measure $\nu$
  on $\A$ such that $\nu_T$ is Haar measure on $G_T$ for any tassel
  $T$ and any choice of trivialization of the endpoints, and $\nu_W=
  \nu_{W^1} \times \cdots \times \nu_{W^k}$ for any web $W$ consisting
  of tassels $W^1, \ldots, W^k$.
\end{theorem}

{\em Proof.\/} We first must prove that $\nu_T$, for a tassel $T$
based at a point $p$, is independent of the choice of trivialization.
A change in the trivialization would effectively replace the holonomy
$g \in G$ of a given connection along $c_i$ by $h_p g_i h_i$, where
$h_p$ and $h_i$ are elements of $G$ expressing the change of
trivialization at the point $p$ and the right endpoint of $c_i$
respectively.  Thus $G_T$ gets sent to $\vec{h}_l
G_T \vec{h}_r$,
where $\vec{h}_l= (h_p, \ldots, h_p)$ and $\vec{h}_r=(h_1, \ldots,
h_n)$, and $h_i=h_j$ if $c_i$ and $c_j$ have the same right endpoint.

Now consider any point $q$ in $\range(T)$.  The set of $t$ such that
all points in $\range(T)$ with parameter value $t$ are regular is open
and dense, so there are such $t<t(q)$ and $t>t(q)$ arbitrarily close
to $t(q)$.  For $t<t(q)$ sufficiently close, every curve that goes
through $q$ goes through exactly one of the regular points with
parameter value $t$, and none shares a regular point with a curve that
does not go through $q$.  Thus the group generated by their types
includes points in $G^n$ with a $g$ in the $i$th entry if $c_i$ goes
through $q$ and $ 1$ if it does not.  Likewise, taking $t>t(q)$ and
small enough, we can find in the group generated by the types elements
in $G^n$ with a $g$ in the $i$th entry if $c_i$ goes through $q$ and
does not end there, and a $1$ otherwise.  Putting these together, we
see that $G_T$ contains every element of $G^n$ with a $g$ in the $i$th
entry if $c_i$ has $q$ as its endpoint and a $1$ otherwise.  Thus
$\vec{h}_r$ is an element of $G_T$. Likewise $\vec{h}_l \in G_T$.  But
since Haar measure on a Lie group is invariant under left and right
multiplication, it gets sent to itself under the map $x \mapsto
\vec{h}_l x \vec{h}_r$.  Thus the assignment of measures to tassels,
and hence to webs, is independent of the choice of trivialization, and
therefore well-defined.

\begin{figure}[hbt]
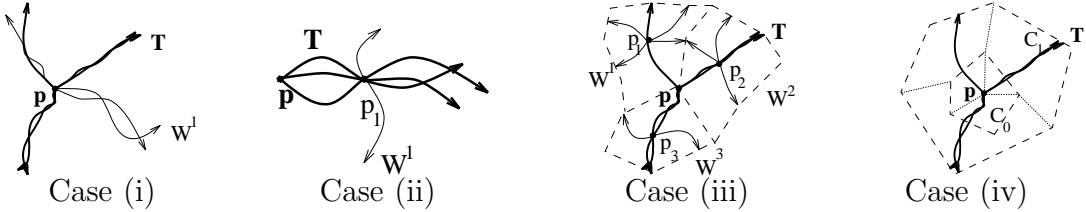

\begin{centering}
\begin{tabular}{cccc}
\pic{case1}{65}{-25}{0}{0} \hspace{8pt}& \hspace{6pt}
\pic{case2}{55}{-25}{0}{0} \hspace{6pt} & \hspace{8pt}
\pic{case3}{65}{-25}{0}{0} \hspace{8pt} & \hspace{8pt}
\pic{case4}{65}{-25}{0}{0} \\ Case (i) & Case (ii) & Case (iii) &
Case (iv)
\end{tabular}
\end{centering}
\caption{Four cases of writing a tassel in terms of a web} \label{fg:cases}
\end{figure}

Now, to check condition (a) of Theorem \ref{thm:gen}, first consider a
tassel $T$ based on $p$, and let $W=\{W^1, \ldots, W^k\}$ be a web on
which $T$ depends. We will show that $\nu_T=(p_{TW})_*\nu_W$ in four
cases.  These are illustrated in Figure \ref{fg:cases}, where the
curves of $T$ are represented in bold and the curves of $W$ are
represented in medium weight.

\begin{romanlist}
\item {\em $W$ consists of a single tassel $W^1$ based at $p_1=p$.}
Since each curve in $W^1$ has $p$ as a left endpoint, and each curve
in $T$ contains $p$ only as its left endpoint, every curve in $T$ is a
curve in $W^1$.  Thus writing $G_T \subseteq G^n$ and $G_{W^1}
\subseteq G^{n_1}$ in the standard way, $p_{TW^1}\maps G_{W^1} \to
G_T$ sends $(g_1, \ldots, g_{n_1})$ to $(g_{k_1}, \ldots, g_{k_n})$
for some integers $1 \leq k_j \leq n_1$.  In particular, $p_{TW^1}$ is
an onto homomorphism from $G_{W^1}$ to $G_T$.  The image of Haar
measure under an onto homomorphism is Haar measure again, so $\nu_T=
(p_{TW^1})_\ast \nu_{W^1}= (p_{TW})_\ast \nu_W$.

\item {\em $W$ consists of a single tassel $W^1$, and $p_1 \neq p$.}
    Let $W'$ be the set of curves in $W^1$ which contain $p$, and
    $W''$ be the set of curves which do not.  Then since every curve
    in $T$ contains $p$ only as its left endpoint, every curve in $T$
    can be written either as $c_i^{-1}$ for $c_i \in W'$ or as
    $c_jc_i^{-1}$ for $c_i \in W'$, $c_j \in W''$.  Clearly
    $(p_{TW^1})_* \nu_{W^1}$ is a probability measure on $G_T$, so it
    suffices to show that $(p_{TW^1})_* \nu_{W^1}$ is invariant under
    right multiplication by elements of $G_T$.  Equivalently, since
    $G_T$ is generated by the types in $T$, we must show that
    $(p_{TW^1})_* \nu_{W^1}$ is invariant under right multiplication
    by any element of $\tau$, for $\tau$ a type of $T$.  For this,
    choose a point $q \in \range(T)$ with $\tau_q =\tau$ and with
    $t(q)$ small enough that $q$ is not on any curve of $W''$.  We can
    identify $\tau$ with $G$ by the canonical isomorphism, and hence
    with $\tau'$, the type of $q$ in $W'$.  Then we have, for $g \in
    \tau$ and $h \in G_{W^1}$
$$p_{TW^1}(h)g=p_{TW^1}(g^{-1}h)$$ so
$$(p_{TW^1})_*(\nu_{W^1})g= (p_{TW^1})_*(g^{-1}\nu_{W^1} ) =
(p_{TW^1})_*(\nu_{W^1}) .$$

\item {\em Each $W^i$ contains a curve containing $p$ that is a
    subcurve of a curve in $T$.}  If there is a $j$ with $p_j=p$, then
    since in a web no tassel is based on a point of intersection with
    other tassels, there is only one $W^j$ in $W$, and we are in
    situation (i).  So assume $p_j \neq p$ for all $j$.

Suppose $c$ in $T$ is product of curves including one in $W^j$ and one
in $W^i$, the one in $W^j$ being the one which contains $p$.  Then
$p_j$ and $p_i$ lie on $c$.  Since $W^i$ contains a subcurve of some
curve $c'$ in $T$ containing $p_j$ and $p$, we have that $c$
intersects $c'$ at a point $p_j \neq p$.  Thus they intersect
infinitely many times, arbitrarily near $p$, and thus the curves in
$W^i$ and $W^j$ intersect infinitely many times.  Since this is
impossible in a web, we conclude that each curve $c$ in $T$ depends on
one $W^i$.  Further, a curve depending on $W^j$ cannot intersect a
curve depending on $W^i$ for $j \neq i$ (except at $p$), because then
they would intersect infinitely many times, arbitrarily near $p$, and
so would curves in $W^j$ and $W^i$ .

Thus the $W^j$'s separate $T$ into families of curves $T^1, \ldots
,T^k$ that intersect only at $p$, with each $T^j$ depending on $W^j$.
But then each $T^j$ is a tassel based at $p$.  Any type of $T$ is a
type of some $T^j$, and commutes with all types of all other $T^i$.
Thus $G_T= G_{T^1} \times \cdots \times G_{T^k}$, and $\nu_T=
\nu_{T^1} \times \cdots \times \nu_{T^k}$ (Of course, $\{T^j\}$ is not
a web, because they all intersect at their bases).  Thus $p_{TW}$ can
be written as a product of maps $p_{T^jW^j}$, so it suffices to show
that $(p_{T^jW^j})_*\nu_{W^j}= \nu_{T^j}$ in order to conclude that
$\nu_{T^1} \times \cdots \times \nu_{T^k}= (p_{TW})_* \nu_W$. But
$T^j$ is a single tassel depending on the single tassel $W^j$, so by
(i) and (ii) we have $\nu_{T^j}= p_{T^jW^j} \nu_{W^j}$.

\item {\em $W$ is arbitrary}.  Let $W_0$ be the set of all $W^j$ which
contain a curve containing $p$ that is a subcurve of a curve in $T$,
and let $W_1$ be the set of all other $W^j$.  Let $C_0$ be $T$
restricted to $\range(W_0)$, and $C_1$ be $T$ restricted to
$\range(W_1)$.  Since every curve in $C_0$ contains $p$, every curve
in $C_1$ does not, and every curve in $T$ contains $p$ exactly once,
it follows that every curve in $T$ is either a curve in $C_0$ or a
product of a curve in $C_1$ and a curve in $C_0$.  Now $C_1$ depends
on $W_1$, so by trivializing $P$ over a neighborhood containing $T$
the set $p_{C_1W_1}(\A_{W_1})$ may be thought of a subset of a product
of copies of $G$.  Since this subset consists of products of types in
$T$, it is contained in $G_T$, so $(p_{C_1W_1})_\ast \nu_{W_1}$ may be
viewed as a probability measure $\mu$ on $G_T$.  Since in this
interpretation $p_{TW_0 \cup W_1}(x_0,x_1) = p_{C_1W_1}(x_1)\cdot
p_{C_0W_0}(x_0)$, the measure $(p_{TW})_*\nu_W$ is the convolution of
$\mu$ and $(p_{C_0W_0})_*\nu_{W_0}$.

Now $C_0$ is a tassel based on $p$, and $G_{C_0}=G_T$, because every
type occurring in $T$ occurs arbitrarily close to $p$, and hence in
$C_0$.  So by (iii) $(p_{C_0W_0})_*\nu_{W_0}= \nu_T$.  But it is well
known that the convolution of a probability measure on a group with
Haar measure is again Haar measure, so $(p_{TW})_*\nu_W=\nu_T$.
\end{romanlist}

To finish checking condition (a) of Theorem \ref{thm:gen}, we suppose
that $W = \{W^1,\dots, W^l\}$ is a web depending on the web $X =
\{X^1, \dots, X^k\}$, and show that $(p_{WX})_\ast \nu_X = \nu_W$.

To see this, note that any $X^i$ can be divided into equivalence
classes $X^i_1, \ldots, X^i_{n_i}$ of curves which are parallel at
$p_i$, and that curves from different equivalence classes do not
intersect except at $p_i$ (this is essentially the argument in point
(iv) above).  Thus every type of $X^i$ is a type of some $X^i_j$, and
commutes with all types of any other $X^i_{j'}$.  In particular
$G_{X^i}= G_{X^i_1} \times \cdots \times G_{X^i_{n_i}}$ and $\nu_{X}=
\nu_{X^1_1} \times \cdots \times \nu_{X^1_{n_1}} \times \nu_{X^2_1}
\times \cdots \times \nu_{X^k_{n_k}}$.

By (i-iv), it suffices to show that $(p_{WX})_* \nu_X$ assigns an
independent measure to each $W^m$, and by the above it suffices to
show that curves in different $W^m$'s do not depend under $p_{WX}$ on
curves in the same $X^i_j$.  This is clear, because if they did then
$p_i$ would be in the range of both of these $W^m$'s, but no neighborhood of it
could be an interval because the curves from the two different tassels
would be parallel at $p_i$.

Condition (b) of Theorem \ref{thm:gen} is immediate.  Each $\nu_W$ is
a probability measure, so as a linear functional it has norm $1$.
\qed

We call this generalized measure $\nu$ the {\it uniform generalized
measure}. This generalized measure has a number of important
properties.  First notice that the group $\Aut(P)$ of automorphisms of
the bundle $P$ acts on the space $\A$, and thus acts as automorphisms
of the C*-algebra $\Fun(\A)$ via
\[ (gF)(A) = F(g^{-1}A) . \] As a consequence it acts dually on the
space $\Fun(\A)^\ast$ of generalized measures on $\A$.  We shall show
that $\nu$ is invariant under this action.  Moreover, since a
generalized measure $\mu$ on $\A$ is equivalent to a measure on
$\overline\A$, it is natural to speak of $\mu$ being a {\it
probability measure} if $\mu(1) = 1$ and for all $F \in \Fun(\A)$, $F
\ge 0$ implies $\mu(F) \ge 0$.  Borrowing some terminology from
C*-algebra theory, we also say that a probability measure $\mu$ is
{\it faithful} if $F \ge 0$ and $\mu(F) = 0$ imply $F = 0$ for all $F
\in \Fun(\A)$.

\begin{corollary}\label{cr:nu}
The uniform generalized measure $\nu$ is a faithful probability
measure, invariant under the action of $\Aut(P)$.
\end{corollary}

{\em Proof.\/} To see that $\nu$ is a faithful probability measure it
suffices to check that $\nu_W$ is a faithful probability measure for
each web $W$.  For this, in turn, it suffices to check it for a
tassel, and Haar measure is clearly a faithful probability measure.
To see that $\nu$ is invariant, note that every step in its
construction was manifestly invariant except the choice of
trivialization, and we showed that $\nu$ was independent of that. \qed

\section{Spin Networks}

Since $\nu$ is a faithful probability measure, we may define $L^2(\A)$
as the Hilbert space completion of the space $\Fun(\A)$ with the inner
product $\bracket{f,g}= \nu(\overline f g)$.  Equivalently, we could
set $L^2(\A) = L^2(\overline A,d\nu)$.  Since $\nu$ is invariant under
$\Aut(P)$, there is a unitary representation of $\Aut(P)$ on
$L^2(\A)$, and thus a unitary representation of the subgroup $\G
\subseteq \Aut(P)$ consisting of gauge transformations.  We define
$L^2(\A/\G)$ to be the closed subspace consisting of vectors in
$L^2(\A)$ invariant under the action of $\G$.  In this section we
describe an explicit set of functions on $\A$ spanning the Hilbert
space $L^2(\A/\G)$; by analogy with the analytic case \cite{Baez5} we
call these `spin networks'.

Given any family $C$, let $L^2(\A_C)$ be the Hilbert space of
square-integrable functions on $\A_C$ with respect to the measure
$\nu_C$.  Of course, the map $f\mapsto f \circ p_C$ from $\Fun(\A_C)$
to $\Fun_0(\A)$ extends to an isometry of $L^2(\A_C)$ into $L^2(\A)$,
and the union of the images of these isometries over all families $C$
is dense in $L^2(\A)$.  In fact if $C$ depends on $D$ then the
embedding of $L^2(\A_C)$ in $L^2(\A)$ factors through that of
$L^2(\A_D)$, so the union of the images of $L^2(\A_W)$ over all webs
$W$ is also dense.  In keeping with the philosophy of this paper, one
can try to understand $L^2(\A)$ by understanding $L^2(\A_W)$ for all
webs $W$.

If $W=\{W^1, \ldots, W^k\}$ is a web, then $L^2(\A_W)$ is fairly
simple to describe.  Fixing a trivialization of $P$ over the endpoints
of the curves, $L^2(\A_W) \iso L^2(G_{W^1}) \tnsr \cdots \tnsr
L^2(G_{W^k})$.  Note however that this isomorphism changes when we
change the trivialization.  Understanding how it changes is a part of
what we need to describe the gauge invariant subspace.

For each family $C$, the group $\G$ acts on $\A_C$.  The quotient of
$\G$ by the subgroup which acts trivially on $\A_C$ is a
finite-dimensional Lie group $\G_C$, which is actually the product
over all endpoints $q$ of curves in $C$ of the groups $\G_q$ of gauge
transformations of the fibers $P_q$.  Fixing a trivialization of $P_q$
gives an isomorphism between $\G_q$ and $G$, so we can think of $\G_C$
as a product of copies of $G$.  The action of $\G$ on $\A_C$ gives a
unitary representation on $L^2(\A_C)$, and when $C$ depends on $D$ the
natural embedding $L^2(\A_C) \hookrightarrow L^2(\A_D)$ is an
intertwining operator.  Let $L^2(\A_C/\G_C)$ be the subspace of
$\G_C$-invariant vectors in $L^2(\A_C)$.  As before, $L^2(\A_C/\G_C)$
embeds into $L^2(\A_D/\G_D)$ if $C$ depends on $D$, they both embed
into $L^2(\A/\G)$ in a consistent fashion, and the image of all such
embeddings is dense in $L^2(\A/\G)$ as $C$ ranges over all families,
or all webs.  We will construct an orthonormal basis of
$L^2(\A_W/\G_W)$ for each web $W$.  The resulting set of vectors will
thus give a set spanning $L^2(\A/\G)$.

To do this, we need an understanding of the action of $\G_W$ on
$L^2(\A_W)$.  We begin by considering the action of $\G_T$ on $\A_T$
when $T$ is a tassel.  If $T$ is a tassel based at $p$, then $\G_p$ is
the group $G$, with action inherited from the left action of $G_T$ on
$L^2(G_T)$ by the map $g \mapsto (g, \ldots , g) \in G_T$.  If $q$ is
any right endpoint of curves in $T$, then $\G_q$ is the group $G$,
with action inherited from the right action of $G_T$ on $L^2(G_T)$ by
the map $g \mapsto (g_1, \ldots, g_n) \in G_T$, where $g_i$ equals $g$
in every entry corresponding to a curve with endpoint $q$, and equals
$1$ otherwise (see the proof of Theorem \ref{thm:uni}).

More precisely, the Peter-Weyl theorem states that $L^2(G_T)$ as a
left and right $G_T$-module decomposes as
$$\bigoplus_{\lambda \in \Lambda_{G_T}} R_\lambda \tnsr
R_\lambda^\dagger,$$ where $\Lambda_{G_T}$ is the set of all
isomorphism classes of irreducible unitary representations of $G_T$,
$R_\lambda$ is an element of the isomorphism class $\lambda$ as a left
representation, and $R_\lambda^\dagger$ is the dual space of
$R_\lambda$, as a right representation.  If $p$ is the base of $T$,
and $H_p$ is the subgroup of $G_T$ consisting of all $(g, \ldots, g)
\in G^n$, then the action of $\G_p$ on $L^2(\A_T) \iso
\bigoplus_{\lambda \in \Lambda_{G_T}} R_\lambda \tnsr
R_\lambda^\dagger $ is the left action of $H_p \subset G_T$.  Likewise
if $H_q$ is the subgroup of $G_T$ consisting of all $(g_1,\dots, g_n)
\in G^n$ with $g_i = g$ if the $i$th curve in $T$ has $q$ as its right
endpoint, and $g_i = 1$ otherwise, then the action of $\G_q$ is the
right action of $H_q \in G_T$ on $\bigoplus_{\lambda \in
\Lambda_{G_T}} R_\lambda \tnsr R_\lambda^\dagger$.

If $W$ is a web, we can write
$$L^2(\A_W) \iso \bigotimes_{j=1}^k \bigoplus_{\lambda_j \in
  \Lambda_j} R_{\lambda_j} \tnsr R_{\lambda_j}^\dagger,$$ where
  $\Lambda_j$ is shorthand for $\Lambda_{G_W^j}$.  The action of the
  gauge group will be the same, except if a point $q$ is the right
  endpoint of more than one tassel, in which case it is the right
  endpoint of two tassels, say $W^j$ and $W^i$.  In this case
  $\G_q\iso G$ acts on $L^2(G_{W^j}) \tnsr L^2(G_{W^i})$ by the tensor
  product of the actions on each individually.  Invariant vectors
  under this action come from invariant elements of
  $R_{\lambda_j}^\dagger \tnsr R_{\lambda_i}^\dagger$ for some choice
  of $\lambda_j$ and $\lambda_i$.  Since the actions of different
  groups $\G_q$ commute, we can decompose each $R_{\lambda_j}^\dagger$
  into an orthogonal direct sum of tensor products, over every $q$ an
  endpoint for $W^j$, of irreducible unitary right representations of
  $\G_q$.

To construct actual $\G_W$-invariant elements of $L^2(\A_W)$, recall
that the Peter-Weyl isomorphism is given by sending the element $v
\tnsr w \in R_\lambda\tnsr R_\lambda^\dagger$ to the function $f(g)=
(w,gv)$ for $g \in G_C$, $(\cdot,\cdot)$ being the usual pairing of a
vector space with its dual, but multiplied by the square root of
$\dim(R)$ to make the isomorphism unitary.  So choose a representation
$\lambda_j \in \Lambda_j$ for $1 \leq j \leq k$, choose an
$H_{p_j}$-invariant vector $\vec{v}_j$ in $R_{\lambda_j}$, choose a
term in the direct sum decomposition of each $R_{\lambda_j}^\dagger$,
such that the representation assigned to each endpoint $q$ which is
only an endpoint for $W^j$ is assigned the trivial representation and
the representations assigned to a $q$ which is an endpoint for $W^i$
and $W^j$ respectively are dual representations.  Also choose a vector
$\vec{w}_q$ in the trivial representation chosen for each $q$ bounding
one tassel, and an invariant element $\vec{w}_q$ of the representation
$V \tnsr V^*$ chosen for each $q$ bounding two tassels.  Notice that
$\bigotimes_{j=1}^k g_j\vec{v}_j$ for $g_j \in G_{W^j}$ is an element
of $\bigotimes_{j=1}^k R_{\lambda_j}$, and that $\bigotimes_{q}
\vec{w}_q$, the product being over all endpoints $q$, is after
reordering appropriately an element of $\bigotimes_{j=1}^k
R_{\lambda_j}^\dagger$, and thus they can be paired (by the rescaled
pairing) to get a number, which we call
$f_{\{\vec{v}_j\},\{\vec{w}_q\}}(g_1, \ldots, g_k)$.  We call such a
function a {\em spin network\/}.
\begin{theorem} \label{thm:spin}\mbox{}
\begin{alphalist}
\item $f_{\{\vec{v}_j\},\{\vec{w}_q\}}$ is invariant under gauge
  transformations, and therefore in particular does not depend on the
  choice of trivialization at endpoints.
\item the function $f_{\{\vec{v}_j\},\{\vec{w}_q\}}$ is in
  $L^2(\A_W/\G_W)$, and furthermore
$$\bracket{f_{\{\vec{v}_j\},\{\vec{w}_q\}},f_{\{\vec{v}'_j\},\{\vec{w}'_q\}}}=
   \prod_{j,q}\bracket{\vec{v}_j,\vec{v}'_j}\bracket{\vec{w}_q,\vec{w}'_q}.$$
\item Choosing $\vec{v}_j$ from an orthonormal basis of the subspace
  of $R_{\lambda_j}$ of $\G_{p_j}$-invariant vectors, and choosing
  fixed unit vectors $w_p$ for each term in the direct sum
  decomposition of $R_{\lambda_j}^\dagger$, we get an orthonormal
  basis for $L^2(\A_W/\G_W)$.
\end{alphalist}
\end{theorem}

{\em Proof.}

\begin{alphalist}
\item The vectors $\{\vec{v}_j\}$ and $\{\vec{w}_q\}$ are invariant
  under gauge transformations, so $f$ is.
\item That $f$ is in $L^2(\A_W/\G_W)$ follows from the previous point
  and the formula for the inner product, which is simply the statement
  that the Peter-Weyl isomorphism is a Hilbert space isomorphism.
\item They are clearly orthonormal, and they certainly span the space
  of spin networks.  But by the Peter-Weyl theorem, every
  $\G_W$-invariant element of $L^2(\A_W)$ is spanned by those of the
  form $\bigotimes_{j=1}^k(\vec{w}_j,g_j\vec{v}_j)$, with $\vec{v}_j$
  and $\vec{w}_j$ invariant elements of $R_{\lambda_j}$ and
  $R_{\lambda_j}^\dagger$ for some $\lambda_j$.  Since such
  $\bigotimes_j \vec{w}_j$ are certainly spanned by all the tensor
  products $\bigotimes_p \vec{w}_p$ used to construct spin-networks,
  it is clear that the spin networks span $L^2(\A_W/\G_W)$.  \qed
\end{alphalist}

\subsection*{Acknowledgements}

The authors would like to thank Greg Kuperberg and Richard Palais for
information on real-analytic structures, Lisa Raphals for the tassel
terminology, and David Vogan for sharing his knowledge of Lie groups
and Haar measure.  J.\ B.\ would also like to thank Jerzy Lewandowski
for inviting him to speak on this subject at the workshop on Canonical
Quantum Gravity in Warsaw.

\end{document}